\documentstyle[fleqn,12pt]{article}
\begin{document}

\newtheorem{theor}{Theorem}
\newtheorem{lemma}{Lemma}
\newtheorem{example}{Example}
\newtheorem{prop}{Proposition}
\newtheorem{assum}{Assumption}
\newtheorem{coroll}{Corollary}

\newcommand{\re}{{I\!\!R}}
\newcommand{\ponto}{\, \cdot  \, }
\newcommand{\formu}[1]{\begin{eqnarray}\label{#1}}
\newcommand{\formub}{\begin{eqnarray} \nonumber}
\newcommand{\eformu}{\end{eqnarray}}
\newcommand{\igual}{\, = \,}
\newcommand{\intre}{\int_{\re}}
\newcommand{\dl}[1]{ \lambda  (d #1)}
\newcommand{\fatori}{\, ! \, }
\newcommand{\bbb}{{\cal B}}
\newcommand{\ie}{{\it i.e. }}
\newcommand{\proof}{\vskip0.4cm \noindent {\bf Proof: }}
\newcommand{\eproof}{ \mbox{}\hfill$\sqcup\!\!\!\!\sqcap$ \vskip0.4cm \noindent}

\title{The Laplace transform and polynomial approximation in $L^2$
      }
\author{Rodrigo Labouriau 
 \thanks{Department of Mathematics,
               Aarhus University.
        }
      }
\date{March, 2016}
\maketitle

\tableofcontents
 
\section{Introduction}

This short note gives a sufficient condition for having the class of polynomials
dense in the space of square integrable functions with respect to a 
finite measure dominated by the Lebesgue measure in the real line, here denoted by $L^2$.
It is shown that if the Laplace transform of the measure in play is bounded 
in a neighbourhood of the origin, then the moments of all order are finite 
and the class of polynomials is dense in $L^2$.
The existence of the moments of all orders is well known for the case where
the measure is concentrated in the positive real line
(see Feller, 1966), but the result 
concerning the polynomial approximation is original, even thought the proof
is relatively simple. This tool is essential for constructing semiparametric extensions of classic parametric models.

A review on the Laplace transform theory is given in section \ref{lap}.
The main result is proved in section \ref{lspol} and an alternative
stronger condition easier to be verified not involving the calculation of the Laplace transform  
is given in section \ref{decay}.

\section{Basic properties of the Laplace Transform} \label{lap}

In this section we review the basic properties
of the Laplace transform.  
Let $f: \re \longrightarrow [0,\infty )$ be a function
such that for some $s\in\re$ the integral
\formu{lap1}
M(s;f) \igual \intre e^{sx} f(x) \dl{x} 
\eformu
converges.
Here $\lambda$ is a $\sigma$- finite measure on $(\re , \bbb (\re) )$.
The function 
$M (\ponto ; f):\re\longrightarrow [0,\infty ]$
such that for each $s\in\re$, $M(s;f)$ is given
by (\ref{lap1}) is said to be the 
{\it Laplace transform} of $f$.

We now study some properties of  the functions with finite Laplace
transform in a neighborhood of zero.
\begin{prop}
\label{lapprop1}
Let $f:\re\longrightarrow [0,\infty )$ be a continuous function such that
for some $\delta > 0$ and for all $s\in (-\delta , \delta )$
\formu{lap4}
M(s;f) < \infty
\,\, .
\eformu
Then $f$ possesses finite  moments of all orders, \ie
for all $n\in N_0=\{0,1,2, \dots \}$,
$$
\intre x^n f(x) \dl{x} \in \re \,\, .
$$
\end{prop}
\proof
Since for all $s\in (-\delta , \delta )$, $M(s;f)<\infty$,
$e^{\vert sx\vert} \le e^{sx} + e^{-sx}$
and using the series version of the monotone convergence theorem
(see Billingsley 1986  page 214 theorem 16.6
\footnote {The referred theorem states:
       ''If $f_n\ge 0$, then $\int \sum_n f_n d\lambda = \sum_n \int f_n d\lambda$.''.}
) we have
\formu{lap1000}
\nonumber
\infty & > & \intre e^{ \delta x} f(x) \dl{x}
       + \intre e^{-\delta x} f(x) \dl{x}
\\ \nonumber & \ge &
\intre e^{\vert \delta x\vert } f(x) \dl{x}
\\ \nonumber & = &
\intre \left \{ \sum_{k=0}^\infty \frac{\vert\delta x\vert^k}{k\fatori} 
\right \} f(x) \dl{x} 
\\ \nonumber & & 
\mbox{(from theorem 16.6 in Billingsley 1986)}
\\ \nonumber & = &
\sum_{k=0}^\infty \left \{
\intre \frac{\vert\delta x\vert^k}{k\fatori} f(x) \dl{x}
\right \}
\, ,
\eformu
and we conclude that the moments of all orders of $f$ are in $\re$.
\eproof

The notion of Laplace transform can be extended to functions with 
range equal to the whole real line in the following way.
Given a function $f:\re\longrightarrow\re$ we define the positive
and the negative part of $f$ respectively by
\formub
f^+ (\ponto ) \igual f(\ponto ) \chi_{[0,\infty )} \{ f (\ponto ) \}
\mbox { and }
f^- (\ponto ) \igual -f(\ponto ) \chi_{(-\infty , 0]} \{ f (\ponto ) \}
\,\, .
\eformu
Here $\chi_A (\ponto )$ is the indicator function of the set $A$.
We have clearly the decomposition
\formub
f(\ponto ) \igual f^+ (\ponto ) - f^- (\ponto )
\,\, .
\eformu
We define the {\it Laplace transform} of a function $f:\re\longrightarrow\re$
as the function $M (\ponto ; f):\re\longrightarrow [-\infty , \infty ]$ 
given by
\formu{lap5}
M(\ponto ; f) \igual M(\ponto ;f^+) - M(\ponto ;f^-)
\,\, ,
\eformu
provided that at least one of the terms of the right side of (\ref{lap5}) is finite (otherwise the Laplace transform of $f$ is not defined).
\begin{prop}
\label {lapprop2}
 Let $f:\re\longrightarrow\re$ and $\delta >0$ be such that for all 
$s\in [-\delta , \delta]$, $M(s;f)\in\re$.
Then, for all $n\in N$ and all
$s\in (-\delta /2 , \delta /2 )$, we have
\formub
M[s; (\ponto )^n f(\ponto ) ] \in \re \,\, .
\eformu
\end{prop}
\proof
Assume without loss of generality that the function $f$ is nonnegative.
Take an arbitrary $s\in [-\delta /2 , \delta /2]$  and $n\in N$.
By hypothesis, $f$ has finite Laplace transform in a neighborhood of zero;
then, from proposition \ref{lapprop1}, $f$ has finite moments of all
orders, in particular
$$
\intre x^{2n} f(x) \dl{x} \in \re 
\,\, .
$$
Using the Cauchy-Schwartz inequality we obtain
\formu {lap6}
\nonumber
\hspace{-8mm}
\vert M[s ; (\ponto )^n f(\ponto ) ] \vert
& = &
\vert  < e^{(\ponto ) s} , (\ponto )^n f(\ponto ) >_\lambda \vert
\\ \nonumber & = &
\vert  < e^{(\ponto ) s}f^{1/2} (\ponto ) , (\ponto )^n f^{1/2}(\ponto ) >_\lambda \vert
\\ \nonumber
& & \mbox{(Cauchy-Schwartz inequality)}
\\ \nonumber
& \le &
\left \| e^{(\ponto ) s}f^{1/2} (\ponto ) \right \| 
\left \| (\ponto )^n f^{1/2}(\ponto ) \right \| 
\\ \nonumber 
&  = &
\left \{ \intre e^{2sx} f(x) \dl{x} \right \}^{1/2}
\hspace{-1mm}
\left \{ \intre x^{2n} f(x) \dl{x} \right \}^{1/2}
\hspace{-5mm}
< \infty
\eformu
\eproof


\section{Polynomial approximation in $L^2$}\label{lspol}

In this section we give a sufficient condition for having the 
class of polynomials dense in $L^2 (a)$.
Here $a$ is a density with respect to the 
$\sigma$- finite measure $\lambda$ of a positive finite
measure on $(\re , \bbb (\re ) )$ and $L^2 (a) $
is endowed with the usual inner product and
norm denoted by $< \, \cdot \, , \, \cdot \, >_a$
and $\| \,\cdot \, \|_a $ respectively.
The conditions we give will ensure that the measure
$a$ possesses all moments finite, \ie
for all $k\in N$,
$$
\intre x^k a(x) \lambda (dx) \, \in \, \re \, .
$$
In that case we can define the sequence of polynomials
$\{ e_i (\ponto ) \}_{i\in N_0} \subseteq L^2 (a)$ 
as the result of a Gram-Schmidt orthonormalization
process applied to the sequence 
$\{ 1 , (\ponto ) , (\ponto )^2 , ... \}$.
The following theorem gives a sufficient condition for
$\{ e_i (\ponto ) \}$ to be a complete sequence in 
$L^2 (a)$, which implies that the polynomials are dense
in  $L^2 (a)$.

\begin{theor}
\label {thpol}
Let $a:\re\longrightarrow\re$ be a function such that
\formu{pol1}
\forall x\in\re , \, a(x) > 0 ;
\eformu
\formu{pol3}
\exists \delta >0  \mbox{ such that } \,
\forall s\in [-\delta ,\delta],
\,
M(s;a) = \intre e^{sx} a(x) \lambda (dx) < \infty\, .
\eformu
Then the orthonormal sequence $\{e_i (\ponto ) \}_{i\in N_0}$ 
is complete in $L^2 (a)$.
\end{theor}
\proof
First of all we observe that condition (\ref{pol3}) implies
that the measure determined by $a$ possesses finite moments
of all orders (see proposition \ref{lapprop1}).

Let $f:\re \longrightarrow \re$ be a function in  $L^2 (a)$
such that for all $k\in N_0$,
\formu{pol4}
\intre x^k f(x)a(x) \lambda (dx) \igual 0 \, .
\eformu
We prove that $f(\ponto ) = 0 $ $a$-a.e. which implies the
theorem (see Luenberg, 1969, Lemma 1, page 61).

Define for each $k\in N_0$, $t\in [-\delta /2 , \delta /2 ]$ and
$x\in\re$,
\formub
f_k (x) \igual (xt)^k f(x) a(x) \, .
\eformu
We will use a series version of the dominated convergence theorem applied to 
$\{ f_k\}$. In the following we find a Lebesgue integrable function
dominating uniformly (\ie for all $k$) the functions $f_k$,
which will enable us to use the referred theorem.
We have for each $n\in N$, $k\in N_0$, 
$t\in [-\delta /2 , \delta /2 ]$ and $x\in\re$,
\formu{pol5}
\left \vert 
\sum_{k=0}^n
f_k (x) 
\right \vert
& \le &
\sum_{k=0}^n
\vert f_k (x) \vert
\igual
\sum_{k=0}^n
\frac{\vert x t \vert^k }{k\fatori } \vert f(x) \vert a(x)
\\ \nonumber
 & = &
\vert f(x) \vert a(x)\sum_{k=0}^n
\frac{\vert x t \vert^k }{k\fatori }
 \, \le \,
\vert f(x) \vert a(x)\sum_{k=0}^\infty
\frac{\vert x t \vert^k }{k\fatori }
\\
\nonumber
 & = &
\vert f(x) \vert a(x) e^{\vert x t \vert }
 \, \le \,
\vert f(x) \vert a(x) \{ e^{xt} + e^{-xt} \}
\\
\nonumber
 & = &
\left \{ \vert f(x) \vert \sqrt{ a(x)}   \right \}
\left \{ \sqrt{ a(x)} (e^{xt} + e^{-xt}) \right \}
\\
\nonumber
 & = &
g(x)
\, ,
\eformu
where the function $g$ is given, for all $x\in\re$, by
\formu{pol6}
g(x) \igual 
\left \{ \vert f(x) \vert \sqrt{ a(x)}   \right \}
\left \{ \sqrt{ a(x)} (e^{xt} + e^{-xt}) \right \}
\,\, .
\eformu
We prove that the function $g$ is Lebesgue integrable.
For, note that 
$$
\left \| \,
\vert f(\ponto ) \vert \sqrt{ a(\ponto )}  
\right \|_{L^2(\lambda )}^2
= \intre \vert f(x) \vert^2 a(x) \lambda (dx) = 
\| f (\ponto ) \|^2_a < \infty .
$$
Then the first term in the right side of (\ref{pol6}) is in 
$L^2(\lambda )$. 
On the other hand,
$$
\left\| 
\sqrt{ a(\ponto )} e^{(\ponto )t}
\right\|^2_{L^2(\lambda )}
= 
\intre e^{2tx} a(x) \lambda (dx) = M(2t;a) < \infty
$$
and 
$$
\left\| 
\sqrt{ a(\ponto )} e^{-(\ponto )t}
\right\|^2_{L^2(\lambda )}
= 
\intre e^{-2tx} a(x) \lambda (dx) = M(-2t;a) < \infty
\, .
$$
Then the second term in the right side of (\ref{pol6}) is in 
$L^2(\lambda )$.
Using the Cauchy-Schwartz inequality
(see Luenberg, 1969, lemma 1, page 47)
we obtain
\formu{pol6a}
\nonumber
\hspace{-9mm}
\left \vert \intre g(x) \lambda (dx) \right \vert 
& = &
\left\vert 
<
\vert f(\ponto ) \vert \sqrt{ a(\ponto )} 
\, , \, 
\sqrt{ a(\ponto )} \left ( e^{(\ponto )t} + e^{-(\ponto )t} \right )
>_\lambda
\right \vert
\\ \nonumber & \le &
\left \| 
\vert f(\ponto ) \vert \sqrt{ a(\ponto )} \right 
\|_{L^2(\lambda )}
\left \| 
\sqrt{ a(\ponto )} 
\left  ( 
e^{(\ponto )t} + e^{-(\ponto )t} 
\right )
\right \|_{L^2(\lambda )}
\hspace{-2mm}
< 
\infty .
\eformu

Since (\ref{pol5}) holds for each $n\in N$, 
$x\in\re$, $t\in [-\delta /2 , \delta / 2]$ and $g$ 
is Lebesgue integrable we can
use the series version of the dominated convergence
theorem
(see Billingsley, 1986, theorem 16.7 page 214
\footnote{
The theorem states:
''If $\sum_n f_n$ converges almost everywhere and
$\vert \sum_{k=1}^n f_k \vert \le g$ almost everywhere,
where $g$ is integrable, then $\sum_n f_n$ and the $f_n$ 
are integrable, and 
$\int \sum_n f_n d\lambda = \sum_n \int f_n d\lambda$''.
}
)
to obtain
\formu{pol7b}
\nonumber
\hspace{-12mm}
\intre e^{xt} f(x) a(x) \lambda (dx)
& = &
\intre 
\left\{ \sum_{k=0}^\infty \frac{(xt)^k}{k\fatori} f(x) a(x)\right\}
\lambda (dx)
\\ \nonumber
& &
\mbox{(from the series dominated convergence theorem)}
\\ \nonumber 
& = &
\sum_{k=0}^\infty
\left \{
\intre  \frac{(xt)^k}{k\fatori} f(x) a(x) \lambda (dx) 
\right \}
\igual 0
\, \, .
\eformu
We conclude that for all $t\in [-\delta / 2 , \delta /2 ]$,
\formu{pol7}
M[t; f(\ponto ) a (\ponto )] \, = \, 0 \, .
\eformu
We show that (\ref{pol7}) implies that $f(\ponto )= 0$
$a$-a.e. . For,
\formu{pol7a}
\nonumber
\hspace{-6mm}
\| f(\ponto ) \|^2_a 
& = &
\vert < \, f(\ponto ) , 1 >_a \vert
\\ \nonumber 
& = &
\left \vert 
< \, \sqrt{f(\ponto )}  \sqrt{f(\ponto )} 
, e^{ (\ponto )\delta/4} e^{-(\ponto )\delta/4 } >_a
\right \vert
\\ \nonumber 
& = & 
\left\vert 
<\, \sqrt{f(\ponto )}  e^{(\ponto )\delta/4} ,
 \sqrt{f(\ponto )}  e^{- (\ponto )\delta/4} >_a
\right \vert
\\ \nonumber  
& & \mbox{(from the Cauchy-Schwartz inequality)}
\\ \nonumber 
& \le &
\left\|
\sqrt{f(\ponto )}  e^{ (\ponto )\delta/4}
\right \|_a
\,\, 
\left\|
\sqrt{f(\ponto )}  e^{-(\ponto )\delta/4} 
\right\|_a
\\ \nonumber  
& = &
\left \{ 
\intre  f(x) e^{\delta/2 x} a(x) \lambda (dx) 
\right \}^{1/2}
\left \{ 
\intre  f(x) e^{-\delta/2 x} a(x) \lambda (dx) 
\right \}^{1/2}
\\ \nonumber  
& = &
\left \{ M[\delta /2 , f(\ponto ) a (\ponto )] \right \}^{1/2}
\left \{ M[- \delta /2 , f(\ponto ) a(\ponto ) ] \right \}^{1/2}
\\ \nonumber  
& = & \mbox{(from (\ref{pol7}))}
 =  0 \, .
\eformu
\eproof

\newpage


\section{Functions with exponentially decaying tails}\label{decay}

The following proposition gives a sufficient condition
for having the Laplace transform defined in a neighborhood
of zero, which is easy to verify.

\begin{prop}
\label{lapprop3}
Let $f:\re\longrightarrow [0,\infty )$ be a continuous function such that
for some $\delta > 0$ and for all $s\in [-\delta , \delta ]$
\formu{lap2}
\lim_{x\rightarrow +\infty} e^{sx} f(x) 
\igual
\lim_{x\rightarrow -\infty} e^{sx} f(x) 
\igual
0
\ponto
\eformu
Then we have:
\begin{description}
\item{ i)}
For all $s\in  (-\delta , \delta  )$ the Laplace transform
of $f$, $M(s;f)$, is finite.
\item{ii)}
For all $k\in N$,
$$
\lim_{x\rightarrow + \infty} x^k f(x) \igual
\lim_{x\rightarrow - \infty} x^k f(x) \igual 0
$$
\end{description}
\end{prop}
\proof

\noindent
{\it i)}
Take $s\in (-\delta , \delta  )$.
Condition  (\ref{lap2}) implies that there exists $L\in\re_+$ such that
for all $x\in\re \setminus [-L , L ]$, $e^{\delta x} f(x) < 1$ and
$e^{-\delta x} f(x) < 1$. We have then
\formu{lap3}
\nonumber
\hspace{-8mm}
M(s; f) & = & \hspace {-3mm} \intre e^{sx} f(x) \dl{x}
\\ \nonumber
& = & \hspace {-3mm}
\int_{[-L,L] } \hspace{-8mm} e^{sx} f(x) \dl{x}
+
\int_{[L,\infty ) } \hspace {-8mm} e^{sx} f(x) \dl{x}
+
\int_{(-\infty ,-L] } \hspace {-9mm} e^{sx} f(x) \dl{x}
\\ \nonumber
& = & \hspace {-3mm}
\int_{[-L,L] } \hspace {-8mm} e^{sx} f(x) \dl{x}
+ \hspace {-1mm}
\int_{[L,\infty ) } \hspace {-8mm} e^{(s-\delta )x } e^{\delta  x} f(x) \dl{x}
+ \hspace{-1.5mm}
\int_{(-\infty ,-L] } \hspace {-11mm} e^{(\delta -s)x} e^{-\delta  x} f(x) \dl{x}
\\ \nonumber
& \le  & 
\int_{[-L,L] } \hspace {-7mm} e^{sx} f(x) \dl{x}
+
\int_{[L,\infty ) } \hspace {-6mm} e^{(s-\delta )x} \dl{x}
+
\int_{(-\infty ,-L] } \hspace {-8mm} e^{(\delta -s)x} \dl{x}
< \infty \, .
\eformu

\noindent
{\it ii)}
For each $k\in N$
$$
\lim_{x\rightarrow + \infty} x^k f(x) 
\igual
\lim_{x\rightarrow + \infty} \{ e^{-\delta x} x^k \} \{ e^{\delta x} f(x) \}
\igual 0
$$
and 
$$
\lim_{x\rightarrow - \infty} x^k f(x) 
\igual
\lim_{x\rightarrow - \infty} \{ e^{\delta x} x^k \} \{ e^{-\delta x} f(x) \}
\igual 0
$$
\eproof

\end{document}